\documentclass[12pt,british]{amsart}

\usepackage[latin1]{inputenc}
\usepackage[T1]{fontenc}
\usepackage{babel,eucal,url,amssymb,stmaryrd,booktabs,%
enumerate,amscd,paralist}

\theoremstyle{plain}
\newtheorem{theorem}{Theorem}[section]
\newtheorem{proposition}[theorem]{Proposition}
\newtheorem{lemma}[theorem]{Lemma}
\newtheorem{corollary}[theorem]{Corollary}

\theoremstyle{definition}
\newtheorem{definition}[theorem]{Definition}

\theoremstyle{remark}
\newtheorem{remark}[theorem]{Remark}
\newtheorem{example}[theorem]{Example}
\newtheorem*{acknowledgements}{Acknowledgements}

\numberwithin{equation}{section}

\newcommand{\Lie}[1]{\operatorname{\textsl{#1}}}
\newcommand{\lie}[1]{\operatorname{\mathfrak{#1}}}

\newcommand{\so}{\lie{so}}
\newcommand{\SP}{\Lie{Sp}}
\newcommand{\sP}{\lie{sp}}

\newcommand{\SU}{\Lie{SU}}
\newcommand{\su}{\lie{su}}
\newcommand{\Un}{\Lie{U}}
\newcommand{\un}{\lie{u}}

\DeclareMathOperator{\Hol}{Hol}

\DeclareMathOperator{\stab}{stab}

\DeclareMathOperator{\ad}{ad}

\newcommand{\Hodge}{\mathord{\mkern1mu *}}

\newcommand{\LC}{\nabla^g}
\newcommand{\RC}{R^g}
\newcommand{\Nt}{\tilde\nabla}
\newcommand{\Rt}{\tilde R}
\newcommand{\Nh}{\hat\nabla}
\newcommand{\Rh}{\hat R}
\newcommand{\Th}{\hat T}

\newcommand{\gp}{\lie g^\bot}
\newcommand{\gu}{\underline{\lie g}}
\newcommand{\Kg}{\mathcal K(\lie g)}
\newcommand{\curv}{\mathcal K}

\newcommand{\Bone}{\textrm{b}_1}
\newcommand{\Btwo}{\textrm{b}_2}

\newcommand{\norm}[1]{\left\lVert #1\right\rVert}
\newcommand{\real}[1]{\left\llbracket #1 \right\rrbracket}

\begin{document}

\title[Einstein metrics via torsion]{Einstein metrics via intrinsic\\ or
parallel torsion}

\author{Richard Cleyton}
\address[Cleyton]{Department of Mathematics and Computer Science\\
University of Southern Denmark\\
Campusvej 55\\
DK-5230 Odense M\\
Denmark}
\email{cleyton@imada.sdu.dk}

\author{Andrew Swann}
\address[Swann]{Department of Mathematics and Computer Science\\
University of Southern Denmark\\
Campusvej 55\\
DK-5230 Odense M\\
Denmark}
\email{swann@imada.sdu.dk}

\begin{abstract}
  The classification of Riemannian manifolds by the holonomy group of their
  Levi-Civita connection picks out many interesting classes of structures,
  several of which are solutions to the Einstein equations. The
  classification has two parts. The first consists of isolated examples:
  the Riemannian symmetric spaces. The second consists of geometries that
  can occur in continuous families: these include the Calabi-Yau structures
  and Joyce manifolds of string theory. One may ask how one can weaken the
  definitions and still obtain similar classifications. We present two
  closely related suggestions. The classifications for these give isolated
  examples that are isotropy irreducible spaces, and known families that
  are the nearly Kähler manifolds in dimension 6 and Gray's weak holonomy
  G$_2$ structures in dimension 7.
\end{abstract}


\maketitle

\newpage
\tableofcontents

\newpage

\section{Introduction}

Linear connections and the equivalent notion of $G$-structures are
fundamental areas of interest in differential geometry. Their equivalence
goes through the holonomy of the connection and it is a well-known Theorem
by Hano and Ozeki that any subgroup $G$ of the general linear group may be
realised as the holonomy of a connection on some
manifold~\cite{Hano-Ozeki:holonomy}. On the other hand, one has the much
more restrictive class of torsion-free connections and their holonomy. The
classification of the irreducible holonomy representations of torsion-free
connections was completed by Merkulov and
Schwachh\"ofer~\cite{Merkulov-Schwachhoefer:irreducible} and recently
Schwachh\"ofer has given a more algebraic
proof~\cite{Schwachhoefer:Berger}. This article concerns itself exclusively
with Riemannian manifolds and therefore the particular case of Riemannian
holonomy will have special significance.  Riemannian holonomies have
yielded geometric structures such as Calabi-Yau manifolds, Joyce manifolds,
hyperK\"ahler and quaternionic K\"ahler manifolds of great interest in both
mathematics and physics.

If $M^n$ is a Riemannian manifold, the holonomy algebra $\lie g$ acts on
tangent spaces via a representation $V$. This induces actions of $\lie g$
on the spaces of trace-less symmetric two-tensors $S^2_0V$ and algebraic
curvature tensors $\Kg$ with values in $\lie g$.  Apart from the generic
holonomy representation of $\so(n)$ on $\mathbb{R}^n$ and that of $\un(n)$
on $\mathbb{C}^n$ of K\"ahler geometry, all holonomy representation satisfy
that the representations $\Kg$ and $S^2_0V$ have no irreducible summands in
common. Schur's Lemma shows that this happens precisely when the equation
$(\Kg{\otimes} S^2_0V)^G = \{0\}$ is satisfied. As a consequence, a
Riemannian manifold is Einstein as soon as the Lie algebra of its holonomy
group is a proper subalgebra of $\so(n)$ not equal to $\un(n/2)$. Our aim
is to generalise this type of statement to the more general framework of
metric connections and their holonomy.


Among metric connections, the connections that give rise to the Riemannian
holonomy groups may be viewed as precisely those with vanishing intrinsic
torsion.  The torsion of a metric $G$-connection is a tensor taking values
in $\Lambda^2V^*{\otimes}V$, where $V$ represents tangent space as a $G$
module. If $\gp$ denotes the orthogonal complement of the Lie algebra $\lie
g$ of $G$ in $\so(n)$ with respect of the metric then the projection of the
torsion to the image of $V\otimes\gp$ under the anti-symmetrising map
$\delta \colon V^*{\otimes}\so(n) \to \Lambda^2V^*{\otimes} V$ is
independent of the chosen $G$-connection. The tensor thus defined is called
the intrinsic torsion of the associated $G$-structure and encodes its
differential geometric data.

Using the decomposition of $V{\otimes} \gp$ into irreducible $G$ modules
one may classify the geometries induced by a $G$-structure according to
where the intrinsic torsion takes its values, an approach first taken by
Gray and Hervella~\cite{Gray-H:16} for almost Hermitian manifolds and since
used by many others.  

Interesting examples of connections with non-trivial intrinsic torsion have
arisen from Gray's definition of weak holonomy~\cite{Gray:weak}.  These
include the six-dimensional nearly K\"ahler geometry with structure group
$\SU(3)$ and weak holonomy $G_2$ in dimension~$7$. Both geometries give
Einstein metrics.

In this article we take the following approach. We consider all
$G$-structures on Riemannian manifolds with non-trivial intrinsic torsion.
If we consider all metric connections with torsion the Hano-Ozeki Theorem
tells us that any Lie subgroup of $\Lie{O}(n)$ can be realised as the
holonomy group of a metric connection on some Riemannian manifold. We
therefore impose various extra conditions on the $G$-structure and its
intrinsic torsion with the specific aim of obtaining Einstein metrics.  The
results are classifications akin to Berger's in the sense that we arrive at
a discrete family of manifolds made up by the non-symmetric isotropy
irreducible homogeneous spaces and continuous families made up of manifolds
of weak holonomy $G_2$ and nearly K\"ahler six-manifolds.

The contents may be outlined as follows: In section~\ref{sec:2}, we give
some of the basic facts and definitions we will need. In
section~\ref{sec:3}, Theorem~\ref{thm:1} gives three conditions on the
intrinsic torsion and tangent space representation of a $G$-structure that
guarantee a solution to the Einstein equations. The first two of these
conditions have as consequences that the intrinsic torsion must be
skew-symmetric and that $G$ acts irreducibly on tangent spaces. Given that
$G$ acts irreducibly on $V$, we investigate the consequences of the third
condition: $(\Kg\otimes S^2_0V)^G=\{0\}$ in the fourth section. These are
listed in~Theorem~\ref{thm:6}. As a by-product of this investigation, we
obtain an algorithm for computing the space of algebraic curvature tensors
for an arbitrary representation $\lie g\to\so(n)$. Using the results of
section~\ref{sec:4}, Theorem~\ref{thm:3} gives a classification of the
holonomy groups of metric connection with parallel torsion when the
holonomy group acts irreducibly on tangent spaces. Finally, in
Theorem~\ref{thm:8}, we make a similar classification of manifolds with a
$G$-structure with invariant skew-symmetric intrinsic torsion and $G$
acting irreducibly on tangent spaces. The crucial fact used here is that
under these conditions we may conclude that the intrinsic torsion is
parallel.

Gray-Hervella type classifications of spaces with skew-symmetric intrinsic
torsion have been made by A.~Fino~\cite{Fino:torsion}, while the case of
skew-symmetric torsion have been considered by, among others, Friedrich and
Ivanov~\cite{Friedrich:non-integrable, Friedrich-I:skew, Ivanov:Spin7}. 
 
\begin{acknowledgements}
  This paper is based on part of the Ph.D.  thesis~\cite{Cleyton:thesis} of
  the first named author written under the supervision of the second named
  author. We thank Andrew Dancer, Anna Fino and Stefan Ivanov for useful
  comments and remarks. The second author wishes to thank F.~M.~Cabrera and
  M.~D.~Monar for their hospitality and useful discussions during the initial
  stage of this work.
  
  Both authors are members of the \textsc{Edge}, Research Training Network
  \textsc{hprn-ct-\oldstylenums{2000}-\oldstylenums{00101}}, supported by
  The European Human Potential Programme.
\end{acknowledgements}

\section{Riemannian Geometry and Metric Connections}
\label{sec:2}

Let $(M^n,g)$ be a Riemannian manifold. The space $\Lie{O}(M)$ of frames on
the tangent bundle orthonormal with respect to the metric $g$ forms a
principal $\Lie{O}(n)$-bundle, where $\Lie{O}(n)$ acts on the right by
change of basis. Write $\mathbb{R}^n$ for the standard representation of
$\Lie{O}(n)$. Then the tangent bundle can be identified with the associated
vector bundle $\Lie{O}(M) \times_{\Lie{O}(n)} \mathbb{R}^n$ and, similarly,
the bundles of $(p,q)$-tensors on $M$ may be identified with
$\Lie{O}(M)\times_{\Lie{O}(n)} (\mathbb{R}^n)^{(p,q)}$. Sections of bundles
of tensors may be identified with equivariant maps $\sigma \colon
\Lie{O}(M) \to (\mathbb{R}^n)^{(p,q)}$ for which
$\sigma(pg)=g^{-1}\sigma(p)$ for $g\in\Lie{O}(n)$ and these sections or
maps may be decomposed according to the decomposition of
$(\mathbb{R}^n)^{(p,q)}$ into irreducible $\Lie{O}(n)$-submodules. In what
follows we will write $\Lambda^p$ for the $p$-th exterior power of $\mathbb
R^n$ as an $\Lie{O}(n)$-representation.
 
\begin{definition}
  Assume that a Riemannian manifold $(M,g)$ has a structure reduction $P
  \subset \Lie{O}(M)$ to a $G$-structure. Let $V$ be the induced
  representation of $G$ on tangent spaces. When $M$ has such a structure
  reduction we will say that the triple $(M,g,V)$ is a \emph{$G$-manifold}.
  
  If the representation $V$ is irreducible we say that $M$ or, more
  precisely, $(M,g,V)$ is \emph{$G$-irreducible}.
  
  If the structure reduction is given by the holonomy $G$ of a metric
  connection $\Nh$ on $(M,g)$ we will say that $(M,g,V)$ is the
  $G$-manifold \emph{determined by} $\Nh$.
  
  A \emph{$G$-connection} $\Nh$ on a $G$-manifold $(M,g,V)$ is a connection
  for which the holonomy group acts on tangent spaces as a subgroup of
  $G$. A \emph{metric connection} on a Riemannian manifold $(M,g)$ is
  connection with holonomy contained in $\Lie{O}(n)$.
\end{definition}

When $(M,g,V)$ is a $G$-manifold the decompositions of
$\Lie{O}(n)$-re\-pre\-sen\-ta\-tions into $G$-modules let us decompose the
tensors of $M$ even further. In this case we have vector bundle isomorphisms
$T^{(p,q)}M \cong (M,g) \times_{\Lie{O}(n)} (\mathbb{R}^n)^{(p,q)} \cong
P\times_G V^{(p,q)}$ and sections may be thought of as $G$-equivariant maps
$\sigma \colon P \to V^{(p,q)}$. Let $\tau$ be a section of the tensor
bundle $T^{(p,q)}M$ and assume that $\tau$ via the identification with a
$G$-equivariant map $P\to V^{(p,q)}$ actually takes its values in some
submodule $W\subset V^{(p,q)}$. When this is the case we use the convenient
notation $\tau\in W$.

\begin{remark}
  Note that representations are real unless we state otherwise and all Lie
  groups and Lie algebras are subgroups and subalgebras of the orthogonal
  ones. Therefore we are free to identify representations with their duals
  and will do so.
\end{remark}

As an important example, let ${\hat\nabla}$ be a metric connection on
a Riemannian manifold $(M,g)$ and let $(M,g,V)$ be the $G$-manifold
determined by $\Nh$. Since $\Nh$ is a metric connection its holonomy
algebra is a subalgebra of $\lie g\leqslant\so(n)$. The difference
between the metric connection and the Levi-Civita connection $\LC$
therefore forms a tensor taking values in $V {\otimes} \so(n)$:
\begin{equation*}
  \eta:=\LC-\Nh\in V {\otimes} \so(n)
\end{equation*}
The map $\delta\colon V {\otimes} \so(n)\to\Lambda^2V{\otimes} V$ given by
$(\delta\alpha)_XY=\alpha_XY-\alpha_YX$ is an isomorphism mapping
$\eta$ to the torsion $\Th$ of $\Nh$. This justifies using the terms
torsion or torsion tensor for either $\eta$ or $\Th$ interchangeably.

Let ${\Rh}$ be the curvature of $\Nh$:
\[
  {\Rh}_{X,Y} = [{\Nh}_X,{\Nh}_Y]-{\Nh}_{[X,Y]}.
\]  
Then the Ambrose-Singer Theorem~\cite{Ambrose-Singer:holonomy} tells us
that ${\Rh} \in \Lambda^2V {\otimes} \lie g$.

On $M$ the Levi-Civita connection $\LC$ is singled out by requiring
that it is a metric connection which is torsion-free:
\begin{equation*}
  T^g_XY:=\LC_XY-\LC_YX-[X,Y]=0.
\end{equation*}
A consequence of this is that the Riemannian curvature tensor $\RC$
satisfies the first and second Bianchi identities:
\begin{gather*}
  \RC_{X,Y}Z+\RC_{Y,Z}X+\RC_{Z,X}Y=0\\
  \LC_X\RC_{Y,Z}+\LC_Y\RC_{Z,X}+\LC_Z\RC_{X,Y}=0.
\end{gather*}
We write $\Bone$ for the map $\Lambda^2 {\otimes} \so(n) \to \Lambda^3
{\otimes} \Lambda^1$ given by
\begin{equation*}
  (\Bone R)(X,Y,Z)=R(X,Y,Z)+R(Y,Z,X)+R(Z,X,Y)
\end{equation*}
and $\Btwo$ for the map $\Lambda^1 {\otimes} \Lambda^2 {\otimes} \so(n) \to
\Lambda^3 {\otimes} \so(n)$ given by
\begin{equation*}
  (\Btwo R')(X,Y,Z,W)=R'(X,Y,Z,W)+R'(Y,Z,X,W)+R'(Z,X,Y,W).
\end{equation*}
Write $\curv(\so(n))$ for the kernel of $\Bone$ and $\curv'(\so(n))$
for the kernel of $\Btwo$. Then the first and second Bianchi
identities say, respectively, that:
\begin{gather*}
  \RC \in \curv(\so(n)) \quad\text{and}\quad \LC\RC \in
  \curv'(\so(n)).
\end{gather*}
Note that when we identify $\Lambda^2\cong\so(n)$, we have $\mathcal
K(\so(n)) \subset S^2(\so(n))$. 

If the Riemannian holonomy of $g$ is contained in $G$ then, by the
Ambrose-Singer Theorem, $\RC$ takes values in $\Lambda^2V{\otimes}\lie
g$. But $\RC$ also satisfies the first Bianchi identity whence
$R\in\curv(\so(n))$. Thus
\begin{equation*}
  R\in \curv(\so(n))\cap(\Lambda^2V{\otimes}\lie g)=\ker\{\Bone\colon \Lambda^2V{\otimes}\lie g\to\Lambda^3V{\otimes} V\}.
\end{equation*}

\begin{definition}
  The representation 
  \begin{equation}
    \label{eq:11}
    \Kg :=\ker\{\Bone\colon \Lambda^2V{\otimes}\lie g\to\Lambda^3V{\otimes} V\}
  \end{equation}
  is called the \emph{space of algebraic curvature tensors with values
    in $\lie g$}.
\end{definition}

In the presence of a metric connection $\Nh$ with torsion $\eta=\LC-\Nh$ we
may write the Riemannian curvature as
\begin{equation}
  \label{eq:1}
  \RC=\Rh+(\Nh\eta)+(\eta^2),
\end{equation}
where $(\Nh\eta)$ is the anti-symmetrisation of the covariant derivative of
$\eta$ with respect to $\Nh$:
\[(\Nh\eta)_{X,Y}Z:=(\Nh_X\eta)_YZ-(\Nh_Y\eta)_XZ\]
and
\[(\eta^2)_{X,Y}Z:=[\eta_X,\eta_Y]Z-\eta_{\eta_XY-\eta_YX}Z.\]
Since $\RC\in\mathcal K(\so(n))$ and $\LC\RC\in\mathcal K'(\so(n))$ we have
\begin{gather}
  \label{eq:2}
  \Bone\Rh=-\Bone(\Nh\eta)-\Bone(\eta^2)
\intertext{and}\label{eq:4}
  \Btwo\Nh\Rh=-\Btwo\Nh(\Nh\eta)-\Btwo\Nh(\eta^2). 
\end{gather}
We will refer to equations~(\ref{eq:2}) and~(\ref{eq:4}) as,
respectively, the first and second \emph{Bianchi relations}. 

\section{Intrinsic Torsion and Einstein Manifolds}
\label{sec:3}

Let $(M,g,V)$ be the Riemannian $G$-manifold and let $\Nh$ be a
$G$-connection on $(M,g)$. Write $\gp$ for the orthogonal complement of
$\lie g<\so(n)$ and $\eta^{\lie g}$ for the component of the torsion tensor
$\eta$ in $V {\otimes} \lie g$. Then the tensor
\begin{equation*}
  \xi:=\eta-\eta^{\lie g}\in V {\otimes} \gp
\end{equation*}
is independent of the choice of $G$-connection on $(M,g,V)$. Corresponding
to $\xi$ is a connection $\Nt:=\LC-\xi$.

\begin{definition}
  When $(M,g,V)$ is a $G$-manifold we call the $G$-con\-nec\-tion $\Nt$ the
  \emph{minimal connection} and $\xi$ the \emph{intrinsic torsion} of the
  $G$-structure.
\end{definition}

This definition is justified by the fact that $\xi$ may be identified with
the intrinsic torsion of the $G$-structure as defined in the introduction
via the isomorphism $\delta$. The argument above proves

\begin{proposition}
  Let $(M,g,V)$ be a $G$-manifold. Then the minimal connection $\Nt$ is the
  unique $G$-connection $\Nt$ on $M$ such that the torsion tensor
  \begin{equation*}
    \xi=\LC-\Nt\in V{\otimes} \gp.
  \end{equation*}
  Among the $G$-connections on $(M,g,V)$ the connection $\Nt$ has the
  property that it minimises the $L^2$-norm of torsion tensors (on compact
  subsets of $M$).\qed
\end{proposition}

The curvature $\Rt$ of the minimal connection is of course related to the
Riemannian curvature precisely as in equation~\eqref{eq:1} for general
metric connections:
\begin{equation}
  \label{eq:3}
  \RC=\Rt+(\Nt\xi)+(\xi^2)
\end{equation}
and the first Bianchi relation for $\Rt$ is
\begin{equation*}
  \Bone\Rt=-\Bone(\Nt\xi)-\Bone(\xi^2). 
\end{equation*}
Write $\Kg$ for the kernel of the restriction $\Bone\colon\Lambda^2V
{\otimes} \lie g \to \Lambda^3V{\otimes} V$ and $\Kg^{\perp}$ for its
orthogonal complement in $\Lambda^2V {\otimes} \lie g$. Then we may split
$\Rt$ into the respective components $\Rt=\Rt_0+\Rt_1$ and conclude that
since $\Bone$ is injective on $\Kg^{\perp}$, $\Rt_1$ is completely
determined by the intrinsic torsion and its covariant derivative. This
observation forms the central idea in

\begin{theorem}\label{thm:1}
  Let $(M,g,V)$ be a Riemannian $G$-manifold.  Assume that the
  intrinsic torsion takes values in the $G$-submodule $W\subset
  V{\otimes}\gp$. Then $g$ is Einstein when the following conditions are
  satisfied:
  \begin{enumerate}[\upshape (a)]
  \item \hfill $(V{\otimes} W{\otimes} S^2_0V)^G=\{0\}$, \hfill \label{cond:A}
  \item \hfill $(S^2W{\otimes} S^2_0V)^G=\{0\}$, \hfill \label{cond:B}
  \item \hfill $(\Kg{\otimes} S^2_0V)^G=\{0\}$. \hfill \label{cond:C}
  \end{enumerate}
\end{theorem}
\begin{proof}
  Write the curvature tensor of $\LC$ as in formula~\eqref{eq:3}. As
  $\xi\in W$ the tensors $(\Nt\xi)$ and $(\xi^2)$ satisfy
  \begin{gather*}
    (\Nt\xi)\in V{\otimes} W,\\(\xi^2)\in  S^2W.
  \end{gather*}
  Since the component $\Rt_1$ of the curvature of the minimal connection is
  determined by these two tensors we also have $\Rt_1\in V{\otimes} W +
  S^2W$. Therefore the three conditions together, through Schur's Lemma,
  imply that no component of the Riemannian curvature can contribute to the
  trace-free Ricci-tensor, whence $g$ is Einstein.
\end{proof}

The conditions are very strong. Condition~(a) implies, first of
all, that no irreducible $G$-module can occur in the decomposition of both
$V$ and $W$, i.e., that $(V{\otimes} W)^G=\{0\}$.  Therefore $(W {\otimes} (V
{\otimes} S^2V))^G=\{0\}$. By exactness of the sequence of $G$-modules
\begin{equation*}
  0 \longrightarrow \Lambda^3V \longrightarrow
  V {\otimes} \Lambda^2V \longrightarrow S^2V {\otimes} V \longrightarrow
  S^3V \longrightarrow 0
\end{equation*}
this implies that $W \subset \Lambda^3V$.

Condition~(b) implies that $V$ is irreducible. We therefore have the
following
\begin{corollary}\label{cor:1}
  Let $(M,g,V)$ be a Riemannian $G$-manifold for which the intrinsic
  torsion takes values in $W\subset V{\otimes}\gp$. If $V$ and $W$
  satisfy conditions $\textrm{(a)}$ and $\textrm{(b)}$ of
  Theorem~\ref{thm:1} then $(M,g,V)$ is $G$-irreducible and the
  intrinsic torsion is a three-form.\qed
\end{corollary}

In particular, Corollary~\ref{cor:1} ensures that it is not restrictive
to assume that $V$ is irreducible when conditions~(a) and~(b) of
Theorem~\ref{thm:1} are satisfied. Furthermore, when $V$ is an
irreducible representation of a Lie group $G$ it is possible to say
precisely when condition~(c) of Theorem~\ref{thm:1} is satisfied. This
will be the main result of the next section.

\section{Berger Algebras and Algebraic Curvature Tensors}
\label{sec:4}
The question we wish to address now is: Given that $G$ is a Lie group
acting irreducibly on a real vector space $V$, when does the space of
algebraic curvature tensors consist of Einstein tensors only? We will
obtain an answer to this at the end of the section. This will be obtained
using tools of the Riemannian holonomy classification, most notably the
concept of the Berger algebra $\gu$ of the Lie algebra of $\lie g$. Note
that the Lie algebra $\lie g$ of $G$ may act reducibly on $V$ even though
$G$ acts irreducibly. Therefore we need to be able to calculate Berger
algebras and the space of algebraic curvature tensors for reducible as well
as irreducible representations.

\subsection{Facts and Definitions}
First a few words on notation: we write $\so(V)$ for the representation of
$\so(n)$ on $V$ whenever $V$ is a real $n$-dimensional vector space.
\begin{definition}
  Let $\lie g$ be a Lie algebra and $V$ be a faithful representation
  of $\lie g$ as a subalgebra of $\so(V)$. We then call the pair
  $(\lie g,V)$ a \emph{metric representation}. If the representation
  $V$ of $\lie g$ is irreducible we say that $(\lie g,V)$ is
  irreducible.

  The Berger algebra $\gu$ of a metric representation $(\lie g,V)$ is the
  smallest subspace $\lie p$ of $\lie g$ such that $\mathcal K(\lie
  p)=\Kg$.
\end{definition}

We collect some facts about the Berger algebra and the space of algebraic
curvature tensors of a metric representation $(\lie g,V)$. The first two
are elementary consequences of the definitions. For proofs of the latter
three we refer the reader to~\cite{Bryant:status, Bryant:hol-survey,
Schwachhoefer:Berger}.
  
  \newcounter{fact}
  \begin{list}{{\bf{Fact}~\arabic{fact}:}}{\usecounter{fact}
  \setlength{\labelwidth 1.5cm}\setlength{\leftmargin
  2.6cm}\setlength{\labelsep 0.5cm}\setlength{\rightmargin 0.5cm}}
  \item $(\gu,V)$ is a metric representation and ~$\underline{\gu}=\gu$.
  \item The space of algebraic curvature tensor for the representation $V$
  of $\lie g$ satisfies
    \begin{equation*}
      \Kg=S^2(\lie g)\cap\mathcal K(\so(V)).
    \end{equation*}
  \item The Berger algebra $\gu$ of a Lie algebra $\lie g$ is an ideal in
    $\lie g$.
  \item  The Berger algebra satisfies
    \begin{equation*}
      \gu = \{r(\alpha):r\in\Kg,\quad\alpha\in\Lambda^2V\}.
    \end{equation*}
  \item A metric representation $(\lie g,V)$ is a Riemannian holonomy
    representation if and only if ~$\lie g=\gu$.
  \end{list}
  
  For ease of reference we also provide a list of the irreducible
  Riemannian holonomy representations and their associated space of
  algebraic curvature tensors here, see Table~\ref{tab:hol}. For a complex
  representation $U$ the notation $\real{U}$ is used to indicate the real
  representation obtained by restricting scalar multiplication to $\mathbb
  R$. In the table and and hereafter the symbols $L,~H,~E$ are used for the
  standard complex representations of $\un(1)$, $\sP(1)$, $\sP(n)$,
  respectively and $\Lambda^{1,0}$ is used for the standard representations
  of both $\un(n)$ and $\su(n)$. We use the notation $V^d_\lambda$ for the
  irreducible representation of dimension $d$ and highest weight $\lambda$.
  Special names have been given to the spin-representation $\Delta$ of
  $\lie{spin}(7)$, the space of Weyl curvature $W$ and the highest weight
  module $\Sigma^{2,2}_0$ of $S^2(\Lambda^{1,0}) {\otimes}
  S^2(\Lambda^{0,1})$.

\begin{table}[tp]
  \begin{tabular}{@{}ccc@{}}
    \toprule
    \( \lie g \)&\( V \)&\( \Kg \)\\
    \midrule
    \( \so(n) \)&\( \Lambda^1=\mathbb R^n \)&\( W+S^2_0+\mathbb R \)\\
    \( \un(n) \)&\( \real{\Lambda^{1,0}}=\mathbb C^n \)&\(
    \Sigma^{2,2}_0+\Sigma^{1,1}_0+\mathbb R \)\\
    \( \su(n) \)&\( \real{\Lambda^{1,0}}=\mathbb C^n \)&\( \Sigma^{2,2}_0 \)\\
    \( \sP(n)\oplus\sP(1) \)&\( EH=\mathbb H^n \)&\( S^4E+\mathbb R \)\\
    \( \sP(n) \)&\( \real{E}=\mathbb H^n \)&\( S^4E \)\\
    \( \lie{spin}(7) \)&\( \Delta=\mathbb R^8 \)&\( V^{168}_{(0,2,0)} \)\\
    \( \lie g_2 \)&\( V^7=\mathbb R^7 \)&\( V^{77}_{(0,2)} \)\\
    \( \lie g \)&\( \lie p \)&\( \mathbb R \)\\
    \bottomrule
  \end{tabular}
  \caption{The irreducible Riemannian holonomy representations and the
    associated spaces of algebraic curvature tensors. In the last row
    $\lie g$ and $\lie p$ denotes the isotropy algebra and
    representation, respectively, of those irreducible symmetric spaces
    not covered by earlier entries in the table.} 
  \label{tab:hol}
\end{table}

\subsection{Reducible Representations}
We start by considering the following special instance:
\begin{example}\label{ex:1}
  Let $V=V_1\oplus V_2$, where $V_1$ and $V_2$ are both non-trivial.
  Consider the inclusion $\so(V_1) \oplus \so(V_2)\subset\so(V)$. We have
  \begin{equation*} 
    S^2(\so(V_1) \oplus \so(V_2)) = S^2(\so(V_1)) \oplus
  \bigl(\so(V_1){\otimes}\so(V_2) \bigr)\oplus S^2(\so(V_2)).
  \end{equation*} 
  The image of $S^2(\so(V)))$ under $\Bone$ is $\Lambda^4V$ which
  decomposes as 
  \begin{equation*} 
    \Lambda^4V =
    \Lambda^4V_1 \oplus \bigl(\Lambda^3V_1 {\otimes}
    V_2\bigr) \oplus \bigl(\Lambda^2V_1 {\otimes}
    \Lambda^2V_2\bigr) \oplus \bigl(V_1 {\otimes}
    \Lambda^3V_2\bigr) \oplus \Lambda^4V_2.
  \end{equation*}
  Let $e_1,\dots,e_p$ be an orthonormal basis of $V_1$ and
  $f_{1},\dots,f_{q}$ an orthonormal basis of $V_2$. Then the set
  consisting of
  \begin{equation*}  
    (e_i\wedge e_j)\vee(f_k\wedge f_l),
  \end{equation*}
  where $1\leqslant i<j\leqslant p$ and~$1\leqslant k<l\leqslant q$, gives
  a basis of the subspace $\so(V_1){\otimes}\so(V_2) \subset S^2(\so(V))$.
  The images of these tensors under $\Bone$ span the subspace
  $\Lambda^2V_1{\otimes} \Lambda^2V_2$ of $\Lambda^4V$. Therefore
  $\Bone\colon\so(V_1){\otimes}\so(V_2))\to\Lambda^2V_1 {\otimes}
  \Lambda^2V_2$ is an isomorphism. Moreover, $\Bone(S^2(\so(V_1))) =
  \Lambda^4V_1$ and $\Bone(S^2(\so(V_2))) = \Lambda^4V_2$, whence
  \begin{align*} 
  \curv(\so(V_1) \oplus \so(V_2)) &=
  S^2(\so(V_1) \oplus \so(V_2))\cap\curv(\so(V))\\ &= 
  \curv(\so(V_1))\oplus \curv(\so(V_2)).
\end{align*}
\end{example}

More generally, we wish to consider the situation where $(\lie g, V)$ is a
metric representation and $V_1$ and $V_2$ are orthogonal submodules of $V$
such that
\begin{equation*}
  V=V_1\oplus V_2.
\end{equation*}
Let $\pi\colon\lie g\to\so(V)$ be the representation of $\lie g$ on $V$ and
let $\pi_i\colon\lie g\to\so(V_i),~i=1,2$ be the two sub-representations.
Then $\pi=\pi_1+\pi_2$ and, since $\pi$ is faithful, the kernels $\hat{\lie
g}_1=\ker\pi_2$ and $\hat{\lie g}_2=\ker\pi_1$ intersect trivially, so
$\hat{\lie g}_1\oplus\hat{\lie g}_2\unlhd\lie g$. On the other hand, if
$\lie g_i:=\pi_i(\lie g)$ then $\lie g\unlhd\lie g_1\oplus\lie g_2$ via the
inclusion $\pi$.

We consider two extremal cases. First, when $\hat{\lie g}_1=\lie g_1$ and
$\hat{\lie g}_2=\lie g_2$. Then $\lie g = \lie g_1 \oplus \lie g_2
\overset{\pi} \hookrightarrow \so(V_1) \oplus \so(V_2)$, where $\pi(\lie
g_i) \subset \so(V_i)$. The computations of Example~\ref{ex:1} show that in
this case
\begin{equation*}
  \Kg = \left(S^2(\lie g_1)\cap\curv(\so(V_1))\right) \oplus \left(S^2(\lie
  g_2)\cap\curv(\so(V_2))\right) = \curv(\lie g_1)\oplus\curv(\lie g_2).
\end{equation*}
In other words, we have
\begin{lemma}\label{lem:7}
  Let $(\lie g_1,V_1)$ and $(\lie g_2,V_2)$ be metric representations. Then
  $(\lie g_1\oplus\lie g_2, V_1\oplus V_2)$ is a metric representation and
\begin{gather}
  \underline{\lie g_1\oplus\lie g_2} = \underline{\lie g_1} \oplus 
  \underline{\lie g_2},\\
  \curv(\lie g_1\oplus\lie g_2) = \curv(\lie g_1)\oplus\curv(\lie g_2).
\end{gather}\qed
\end{lemma}
The lemma shows that new metric representations may be obtained by making
direct products of Lie algebras and representations and, furthermore, that
the Berger algebra and space of algebraic curvature tensors of these new
metric representation are obtained from those of the summands by direct
product. If a metric representation is obtained in this fashion, we will
use the shorthand notation
\begin{equation}\label{eq:10}
  \bigoplus(\lie g_i,V_i) := \left(\bigoplus\lie g_i, \bigoplus
  V_i\right).
\end{equation}
The second extremal case is when $\hat{\lie g}_1=\{0\}=\hat{\lie g}_2$.

\begin{lemma}\label{lem:6}
  If $(\lie g, V)$ is a metric representation and $V=V_1\oplus V_2$ is an
  orthogonal decomposition of $V$ into submodules $V_1, V_2$ such that both
  the induced representations $\pi_i\colon\lie g\to\so(V_i)$ are faithful
  then
  \begin{equation*}
    \gu=\{0\}\quad\text{and}\quad\Kg=\{0\}.
  \end{equation*}
\end{lemma}
\begin{proof}
  Let $(\lie g,V)$ be a metric representation and $V=V_1\oplus V_2$ is
  an orthogonal decomposition of $V$ where both the
  sub-representations $\pi_i\colon\lie g \to \so(V_i)$ are faithful.
  Then $\lie g$ is included diagonally into a direct sum $\lie g_1
  \oplus \lie g_2$ of two copies $\lie g_i,~i=1,2$ of $\lie g$ with
  $\lie g_i \subset \so(V_i)$. Therefore $S^2(\lie g)\subset S^2(\lie
  g_1)\oplus \lie g_1{\otimes}\lie g_2\oplus S^2(\lie g_2)$, where
  $S^2(\lie g_i)\subset S^2(\so(V_i))$ and $\lie g_1{\otimes}\lie g_1
  \subset \so(V_1){\otimes}\so(V_2)$. As $\Bone(\lie g_1{\otimes}\lie
  g_2) \cong \lie g_1{\otimes}\lie g_2$ is orthogonal to
  $\Bone(S^2(\lie g_1) \oplus S^2(\lie g_2))$ any element $R$ of $\Kg$
  must project to zero in $\so(V_1){\otimes} \so(V_2)$. Assume that
  $\gamma^1,\dots,\gamma^d$ is an orthogonal basis of $\lie g$. Write
  $\gamma^k=\gamma^k_1+\gamma^k_2$, where
  $\gamma^k_i=\pi_i(\gamma^k)$.  Any element $R\in\Kg$ may then be
  written as
  \begin{equation*}
    R=\sum_{k\leqslant l}a_{kl}\gamma^k\vee\gamma^l.
  \end{equation*}
  The projection of $R$ to $\so(V_1){\otimes} \so(V_2)$ is
  \begin{equation*}
    \sum_{k\leqslant l} a_{kl} (\gamma^k_2\vee\gamma^l_1 +
    \gamma^k_1\vee\gamma^l_2) = \sum_{k,l}a_{kl}\gamma^k_1\vee\gamma^l_2
  \end{equation*}
  which is zero only if $R=0$. The Lemma follows.
\end{proof}

\begin{proposition}
  Let $(\lie g, V)$ be a metric representation. Assume that $V=V_1\oplus
  V_2$ is a decomposition of $V$ into orthogonal submodules and let
  $\pi_i:\lie g\to\so(V_i),~i=1,2$ be the induced representations.  Define
  $\hat{\lie g}_1:=\ker\pi_2$ and $\hat{\lie g}_2:=\ker\pi_1$. Then
  $(\hat{\lie g}_1,V_1)$ and $(\hat{\lie g}_2,V_2)$ are metric
  representations such that
  \begin{gather}
    \gu=\hat{\gu}_1\oplus\hat{\gu}_2\\\intertext{and} 
    \Kg=\curv(\hat{\lie g}_1)\oplus\curv(\hat{\lie g}_2).
  \end{gather}
\end{proposition}
\begin{proof}
  When $(\lie g,V)$ is metric representation with an orthogonal
  decomposition of $V$ into submodules $V_1\oplus V_2$ we define $\hat{\lie
  g}_i$, and $\lie g_i$ as before. Let $\tilde{\lie g}\unlhd\lie g$ be the
  orthogonal complement of $\hat{\lie g}_1\oplus\hat{\lie g}_2$ in $\lie
  g$. Then $\lie g_i\cong\hat{\lie g}_i\oplus\tilde{\lie g}_i$, where
  $\tilde{\lie g}_i=\pi_i(\tilde{\lie g})$ and thus 
  \begin{equation*}
    \hat{\lie g}_1\oplus\hat{\lie g}_2\subset\lie g \subset \lie g_1 \oplus \lie g_2 \cong 
    \left(\hat{\lie g}_1 \oplus \tilde{\lie g}_1\right)
    \oplus \left(\hat{\lie g}_2 \oplus \tilde{\lie g}_2\right) 
    \subset\so(V_1) \oplus \so(V_2).
  \end{equation*}
  By Lemma~\ref{lem:7}, 
  \begin{equation*}
    \curv(\hat{\lie g}_1 \oplus \hat{\lie g}_2) =
    \curv(\hat{\lie g}_1) \oplus \curv(\hat{\lie g}_2)
    \quad\text{and}\quad 
    \curv(\lie g_1 \oplus \lie g_2) = \curv(\lie g_1) \oplus
    \curv(\lie g_2),
  \end{equation*}
  whence
  \begin{equation*}
    \curv(\hat{\lie g}_1) \oplus \curv(\hat{\lie g}_2) \subset 
    \Kg\subset\curv(\lie g_1)\oplus
    \curv(\lie g_2)\subset S^2(\lie g_1) \oplus S^2(\lie g_2).
  \end{equation*}
  The final inclusion shows that any algebraic curvature tensor
  $R\in\Kg$ must have trivial projection to $\tilde{\lie
    g}_1{\otimes}\tilde{\lie g}_2\subset S^2(\lie g_1\oplus\lie g_2)$.
  By an argument similar to the one given in the proof of
  Lemma~\ref{lem:6}, any curvature tensor in $\Kg$ must satisfy that
  the component taking values in $\tilde{\lie g}$ vanishes and
  therefore $\Kg = \curv(\hat{\lie g}_1) \oplus \curv(\hat{\lie
    g}_2)$.
\end{proof}

\subsection{Berger Decomposition}

Metric representations are not generally of the form given by
equation~(\ref{eq:10}). An obvious question to ask is therefore: how may
one compute the Berger algebra and the space of algebraic curvature tensors
for an arbitrary metric representation $(\lie g, V)$? The results of the
previous section will allow us to introduce a Berger decomposition of the
metric representation and show that the Berger algebra of the metric
representation may be computed as a direct sum of Berger algebras of the
irreducible summands of the Berger decomposition. The irreducible case is
then dealt with in the next section.

\begin{definition}
  Let $(\lie g,V)$ be a reducible metric representation and let 
  \begin{equation*}
    V=\bigoplus_{i}V_i\tag{$\Hodge$}
  \end{equation*}
  be an orthogonal decomposition of $V$ into irreducible submodules.
  For each $i$ let
  \begin{equation*}
    \hat{V}_i:=\bigoplus_{j\not=i}V_j,
  \end{equation*}
  let $\pi_i$ and $\hat{\pi}_i$ be the projections $\pi_i \colon \lie
  g \to \so(V_i)$ and $\hat{\pi}_i \colon \lie g \to \so(\hat{V}_i)$,
  and let $\lie g_i=\ker\hat{\pi}_i$. Then \emph{the Berger decomposition} of
  $(\lie g, V)$ with respect to the decomposition~($\Hodge$) is
  \begin{equation*}
    B(\lie g, V)=\bigoplus_i(\lie g_i,V_i).
  \end{equation*}
\end{definition}
The definition and the results of the previous section proves
\begin{proposition}\label{prop:3}
  Let $(\lie g, V)$ be a metric representation. Assume that
  \begin{equation*}
    B(\lie g, V)=\bigoplus_i(\lie g_i,V_i).
  \end{equation*}
  is a Berger decomposition of $(\lie g, V)$. Then
  \begin{gather*}
    \gu=\bigoplus_i\gu_i\quad\text{and}\quad\Kg=\bigoplus_i\curv(\lie g_i).
  \end{gather*}\qed
\end{proposition}

\begin{example}
  Let $\lie g$ be a simple Lie algebra and let $V$ be a non-trivial,
  real representation of $\lie g$. Then $(\lie g, V)$ is a metric
  representation. If $V_i\subset V$ is a non-trivial submodule of $V$
  then $\lie g_i$ is non-trivial only if $\hat{V}_i\cong k\mathbb R$, where
  $k$ is a non-negative integer. Therefore $\gu\not=\{0\}$ only if $V$
  is either irreducible or isomorphic to $V'\oplus k\mathbb R$ for some
  irreducible representation $V'$.
\end{example}

Note that the component representations of a Berger decomposition may
themselves be reducible. So to calculate the Berger algebra of an arbitrary
metric representation we might need to invoke Proposition~\ref{prop:3}
several times. However, we also have the following corollary of
Lemma~\ref{lem:6}:
\begin{lemma}\label{lem:5}
  Let $(\lie g, V)$ be a metric representation. If $V\cong kV'$ for $k>1$
  and some representation $V'$ of $\lie g$ then $\gu=\{0\}$. 
\end{lemma}
\begin{proof}
  If $(\lie g,V)$ is a metric representation where $V=kV'$ then any
  orthogonal decomposition $V=V'\oplus(k-1)V'$ satisfies that the
  induced representations $\pi'\colon\lie g\to\so(V')$ and
  $\pi''\colon\lie g\to\so((k-1)V')$ are faithful. Therefore
  Lemma~\ref{lem:6} applies.
\end{proof}

Using Lemma~\ref{lem:5} we may eliminate any reducible component
representations from a Berger decomposition. To see this, consider the case
of $\lie g$ represented on $V = V_1\oplus V_2$ where $V_1$ is an
irreducible representation and $V_2$ is its orthogonal complement.  Write
$\lie g=\hat{\lie g}_1\oplus\hat{\lie g}_2\oplus\tilde{\lie g}$ as above.
Then $\hat{\lie g}_1\oplus\tilde{\lie g}\cong\pi_1(\lie g)$ acts
irreducibly on $V_1$. Thus, $\lie g_1$ must act on $V_1$ as a direct sum of
isomorphic representations if it does not act irreducibly.

The results we have found in the present section, form an algorithm for
finding the Berger algebra of an arbitrary metric representation $(\lie g,
V)$: first decompose $V$ into irreducible submodules $V=\bigoplus V_i$ and
construct its Berger decomposition $B(\lie g, V)=\bigoplus(\lie g_i,V_i)$.
Then the Berger algebra may be computed using Proposition~\ref{prop:3},
Lemma~\ref{lem:5} and Proposition~\ref{prop:5}.

\subsection{Irreducible Representations}

The promised result for determining when condition~(c) of
Theorem~\ref{thm:1} is satisfied, is nearly at hand. In fact, from the
results of the previous section we may conclude

\begin{proposition}\label{prop:5}
  Let $(\lie g,V)$ be an irreducible metric representation. Then $(\lie
  g,V)$ satisfies either
  \begin{enumerate}[\upshape (i)]
    \item $\gu=\{0\}$,\label{item:6}
    \item $\lie g$ is a Riemannian holonomy representation,\label{item:7}
    \item or $\gu=\sP(n)$ and $(\lie g,V)=(\sP(n)+\un(1),\real{EL})$, where
      $E$ and $L$ are the standard complex representations of $\sP(n)$ and
      $\un(1)$.\label{item:8}
  \end{enumerate}
\end{proposition}
\begin{proof}
  Let $(\lie g,V)$ be an irreducible metric representation. Assume
  that $\gu\not=\{0\}$. Write $\lie g=\lie g_1\oplus\lie g_2$, where
  $\lie g_1=\gu$. Then the complexification $V{\otimes}\mathbb C$ falls
  into one of the following cases, depending on the types of
  representation of $\lie g$,~$\lie g_1$ and $\lie g_2$ on $V$. In
  this respect we follow the conventions of Br\"ocker and tom
  Dieck~\cite{Broecker-tom-Dieck:Lie}.
\begin{enumerate}
\item If $V$ is of real type, then either
\begin{enumerate}
\item $V{\otimes}\mathbb C=U_1{\otimes} U_2$, where $U_1$ and $U_2$ are
  irreducible complex representations of real type, or,\label{item:1}
\item $V{\otimes}\mathbb C=U_1{\otimes} U_2$, where $U_1$ and $U_2$ are
  irreducible complex representations of quaternionic
  type.\label{item:2}
\end{enumerate}
\item If $V$ is of complex type, then $V{\otimes}\mathbb C=U_1{\otimes}
  U_2+\overline{U_1{\otimes} U_2}$ where $U_1$ and $U_2$ are irreducible
  complex representations and either $U_1$ or $U_2$ is of complex
  type.\label{item:3}
\item If $V$ is of quaternionic type then $V{\otimes}\mathbb
  C=2U_1{\otimes} U_2$ where $U_1$ and $U_2$ are irreducible complex
  representations and either
\begin{enumerate}
\item $U_1$ is of quaternionic type and $U_2$ is of real type,
  or,\label{item:4}
\item $U_1$ is of real type and $U_2$ is of quaternionic
  type.\label{item:5}
\end{enumerate}
\end{enumerate}
Lemma~\ref{lem:5} ensures that the restriction of the
representation $V$ to $\lie g_1$ is irreducible. This places severe
restrictions on the dimension of $U_2$. To obtain the desired result
all that is needed is essentially book-keeping: In case~(\ref{item:1}),
the dimension of $U_2$ must be one, so $\lie g_2 \leqslant \so(1) =
\{0\}$. In case~(\ref{item:2}) if $q=\dim U_1$ then $\lie g_1=\sP(q/2)$
as this is the only holonomy representation of quaternionic type, and
$\dim U_2 = 2$ whence $\lie g_2 = \{0\}$ or $\sP(1)$. However, the
later possibility contradicts the holonomy classification, so $\lie
g_2=\{0\}$.

In case~(\ref{item:3}), $\dim U_2$ must be one for $V$ to be irreducible
and $\lie g_2$ must then be either $\{0\}$ or $\un(1)$. Then $\lie
g_1$ and $U_1$ must be either $\un(n),~\su(n)$ or $\sP(n)$ acting on
their standard complex representations. However, $\un(n)\oplus\un(1)$
does not act faithfully on $\Lambda^{1,0}{\otimes} L$ and
$\su(n)\oplus\un(1)$ acting on $\Lambda^{1,0}{\otimes} L$ is a holonomy
representation, so the only possibilities are $\lie g_2=\{0\}$ or
$\lie g_1 = \sP(n)$ and $\lie g_2=\un(1)$ with $U_1=E$, $U_2=L$. 

In case~(\ref{item:4}), $\dim U_2$ is one again and $\lie g_2 \leqslant
\so(1) = \{0\}$. In case~(\ref{item:5}), $\dim U_2\geqslant 2$ and this
implies that $V$ is reducible and thus we have a contradiction with
the initial assumption.
\end{proof}

\begin{corollary}\label{cor:2}
  If $V$ is an irreducible representation of a Lie algebra $\lie g$
  then either $\Kg=\{0\}$ or $(\lie g, V)$ is an irreducible holonomy
  representation, or $(\sP(n)\oplus\un(1),\real{EL})$.\qed
\end{corollary}

The assumption of irreducibility in Corollary~\ref{cor:2} is not quite what
we want. If $G$ is a connected Lie group there is no problem as any
irreducible $G$-representation will be an irreducible module of its Lie
algebra $\lie g$. If $G$ is not connected we may have an irreducible $G$
representation $V$ that is reducible as a representation of $\lie g$. But
then its decomposition as a $\lie g$-module is into a direct sum of
isomorphic submodules. This is so since the identity component $G_0$ of $G$
preserves the $\lie g$-irreducible submodules. So if $V_1$ and $V_2$ are
$G_0$-irreducible subrepresentations of $V$ then there is some element of
$G \setminus G_0$ that maps $V_1$ to $V_2$ whereby they are seen to be
isomorphic as $G_0$-representations. Lemma~\ref{lem:5} then yields:
\begin{corollary}\label{cor:3}
  Let $V$ be an $n$-dimensional, irreducible, real representation of a Lie
  group $G$. If the Lie algebra $\lie g$ of $G$ acts reducibly on $V$ then
  $\Kg=\{0\}$.\qed
\end{corollary}
The Corollaries~\ref{cor:2} and~\ref{cor:3} and the third column of
Table~\ref{tab:hol} give us:
\begin{theorem}\label{thm:6}
  Let $V$ be an irreducible representation of a Lie group $G$. Then one of
  the following holds
  \begin{enumerate}[\upshape (a)]
  \item $\Kg=\{0\}$,
  \item $G$ acts on $V$ as an irreducible holonomy representation, or,
  \item $G=\SP(n)\Un(1)$ and $V=\real{EL}$.
  \end{enumerate}
  In particular, the space of algebraic curvature tensors consists only of
  Einstein tensors if and only if $\lie g$ is a proper subalgebra of
  $\so(n)$ and $V$ is not the standard representation of $\un(n/2)$.\qed
\end{theorem}
Note for future use, that in case~(c) the Berger algebra is $\sP(n)$.

\section{Parallel Torsion and Einstein Metrics}
\label{sec:5}
In this section we return to the following set-up: Let $(M,g,V)$ be a
$G$-manifold determined by $\Nh$, a metric connection on $M$. Let $\eta \in
V {\otimes} \Lambda^2V$ be its torsion tensor. Assume that the torsion tensor
is parallel with respect to $\Nh$:
\begin{equation*}
  \Nh\eta=0.
\end{equation*}
Note that this implies that $\eta$ is invariant by the holonomy $G$ of
$\Nh$, whence $\eta_X^{\lie g}.\eta=0$ where $.$ denotes the standard
action of $\so(n)$ on $V{\otimes}\so(n)$. We conclude that the
intrinsic torsion $\xi$ of $(M,g,V)$ is invariant by $G$ as well as
parallel with respect to the minimal connection. Therefore the
following definitions are equivalent:

\begin{definition}
  Let $(M,g)$ be Riemannian manifold. We say that \emph{ $(M,g)$ has
    parallel torsion} if it admits a metric connection $\Nh$ for which
  the torsion
  $\eta$ satisfies $\Nh\eta=0$. 
\end{definition}

\begin{definition}
  Let $(M,g,V)$ be a $G$-manifold. We say that \emph{$(M,g,V)$ is a
  parallel $G$-manifold} if the intrinsic torsion is parallel with respect
  to the minimal connection.
\end{definition}

The following result is then an easy consequence of
equation~(\ref{eq:1}) and invariance of the torsion by the holonomy
group.

\begin{theorem}\label{thm:7}
  Let $(M,g,V)$ be a parallel $G$-manifold which furthermore is
  $G$-irreducible. Then $(M,g)$ is Einstein if $(\Kg{\otimes}
  S^2_0V)^G=\{0\}$.\qed
\end{theorem}

\subsection{Ambrose-Singer Manifolds}

A particular instance of parallel $G$-manifolds are those which admit a
connection for which both the curvature and the torsion are parallel.

\begin{definition}
  Let $D$ be a metric connection on a Riemannian manifold $(M,g)$ for which
  the curvature $R^D$ and torsion $T^D$ satisfies
  \begin{equation*}
    DT^D=0,\qquad DR^D=0.
  \end{equation*}
  Then $D$ is called \emph{an Ambrose-Singer connection}. A triple
  $(M,g,D)$ where $(M,g)$ is a Riemannian manifold and $D$ is an
  Ambrose-Singer connection will be called \emph{an Ambrose-Singer
    manifold}.
\end{definition}

\begin{remark}  
  In the literature, an Ambrose-Singer manifold is often called a
  \emph{locally homogeneous manifold}. Note that Ambrose-Singer manifolds
  are \emph{not} locally diffeomorphic to homogeneous spaces.
\end{remark}

To each Ambrose-Singer connection $D$ we may of course associate the
$G$-manifold $(M,g,V)$ given by its holonomy. Thereby we obtain a
parallel $G$-manifold $(M,g,V)$.  An obvious question is: When do we
obtain an Ambrose-Singer manifold from a Riemannian manifold with
parallel torsion?

\begin{lemma}
  Let $(M,g)$ be a Riemannian manifold and let $\Nh$ be a metric connection
  on $M$ with parallel torsion.  If the Riemannian $G$-manifold $(M,g,V)$
  given by the holonomy of $\Nh$ has trivial Berger algebra then
  $(M,g,\Nh)$ is an Ambrose-Singer manifold.
\end{lemma}

\begin{proof}
  Let $(M,g)$ be a Riemannian manifold with a metric connection $\Nh$ for
  which the torsion tensor $\eta$ is parallel. Let $\lie g$ be the Lie
  algebra of the holonomy of $\Nh$. Assume that $\gu=\{0\}$ or,
  equivalently, that $\mathcal{K}(\lie g)=\{0\}$. This implies that the
  curvature $\Rh$ is determined completely by the tensor $(\eta^2)$ through
  the Bianchi relation~\eqref{eq:2}. Therefore both
  \begin{equation*}
    \Nh\eta=0\quad\text{and}\quad\Nh\Rh=0
  \end{equation*}
  hold.
\end{proof}

To an Ambrose-Singer manifold $(M,g,D)$ one may also associate an
\emph{infinitesimal model}. Briefly, this consists in building a Lie
bracket $[ \cdot , \cdot ]_{\lie h}$ on $\lie h:=\lie g\oplus V$ by
defining
\begin{equation}\label{eq:5}
  [A+X,B+Y]_{\lie h}:=\left([A,B]_{\lie g}+\Rh_{X,Y}\right)
                       +\left(AY-BX-\Th_XY\right).
\end{equation}
where $\Th_XY=-\eta_XY+\eta_YX$ is the `usual' torsion of $\Nh$, $A,B \in
\lie g$ and $X,Y \in V$. The Bianchi relations and invariance of $\hat T$
and $\Rh$ by $\lie g$ ensures that this satisfies the Jacobi-identity.
Thus, we obtain a pair of Lie algebras $(\lie g,\lie h)$ with $\lie g
\leqslant \lie h$.

\begin{definition}
  Let $(\lie g,\lie h)$ be a pair of Lie algebras. We say that $(\lie
  g,\lie h)$ is \emph{effective} if $\lie g\leqslant\lie h$ and the
  representation of $\lie g$ on $\lie h/\lie g$ is faithful.
\end{definition}

\begin{remark}
  When $\lie g$ is the holonomy algebra of an Ambrose-Singer
  connection and $\lie h=\lie g\oplus V$ with Lie bracket defined as
  in equation~(\ref{eq:5}) above then the pair of Lie algebras $(\lie
  g,\lie h)$ is effective.
\end{remark}

\begin{definition}
  Let $\lie g$ and $\lie h$ be Lie algebras such that $(\lie g,\lie
  h)$ is effective. Let $H$ be the connected, simply-connected Lie
  group with Lie algebra $\lie h$ and $G$ the connected Lie subgroup
  of $H$ with Lie algebra $\lie g$.  We will say that $(\lie g,\lie
  h)$ is \emph{regular} if $H$ a closed subgroup of $G$.  Similarly we will
  say that \emph{an Ambrose-Singer manifold is regular} if the pair of Lie
  algebras obtained from its infinitesimal model is regular.
\end{definition}

In~\cite{Tricerri:Local-homogeneous}, Tricerri proved the following
Theorem.

\begin{theorem}\label{thm:2}
  An Ambrose-Singer manifold is locally isometric to a homogeneous space if
  and only it is regular.\qed
\end{theorem}

\begin{corollary}
  Let $(M,g,V)$ be a $G$-manifold where $G$ is the holonomy of a metric
  connection $\Nh$ on $M$ with parallel torsion. Assume that the Berger
  algebra of $\lie g$ is trivial. Then $(M,g,\Nh)$ is an Ambrose-Singer
  manifold and $(M,g)$ is locally isometric to a homogeneous space if and
  only if $(M,g,\Nh)$ is regular.\qed
\end{corollary}

Note that if a Lie algebra $\lie h$ has a reductive decomposition
$\lie h=\lie g\oplus V$ where $V$ is an irreducible and faithful
representation of $\lie g$ then $\overline{G}$ must be either $H$ or
$G$. By continuity $\overline{G}$ preserves the reductive
decomposition $\lie g\oplus V$. So assuming that $\overline{G}=H$
leeds to the conclusion that the action of $\ad_{\lie h}$ preserves
the subspaces of the reductive splitting. In particular, $[\lie
g,V]\subset V$ and $[V,\lie g]\subset \lie g$. This implies that
$\lie g$ acts trivially on $V$ and therefore establishes a
contradiction. So we have

\begin{proposition}\label{fn:1}
  Let $(M,g,V)$ be a $G$-irreducible manifold where $G$ is the holonomy of
  a metric connection for which the torsion is parallel. Assume that the
  Berger algebra of $\lie g$ is trivial. Then $(M,g)$ is locally isometric
  to an isotropy irreducible homogeneous space $H/G$.\qed
\end{proposition}

\subsection{Classification}

Let $(M,g)$ be a Riemannian manifold, parallel with respect to some metric
connection $\Nh$.  Assume that the holonomy algebra $\lie g$ of $\Nh$ acts
irreducibly on the tangent spaces $V$. If the torsion is assumed to be
non-trivial, this immediately places heavy restrictions on the pair $(\lie
g,V)$ since $\lie g$ must leave some tensor in $V {\otimes} \so(n)$
invariant. If we write $\so(n)=\lie g \oplus \gp$ then we have one of two
possibilities: Either $\lie g\cong V$ and $\lie g$ is a simple Lie algebra
or $\gp$ contains a submodule isomorphic to $V$.  The following lemma is
obtained by inspection of representations.

\begin{lemma}\label{lem:1}
  Let $V\cong\mathbb{R}^n$ be an irreducible representation of a Lie
  algebra $\lie g\leqslant\so(n)$ and let $\gp$ be the orthogonal
  complement of $\lie g$ in $\so(n)$.
  \begin{enumerate}[\upshape (i)]
  \item Assume that $\left(V {\otimes} \lie g\right)^G \not= \{0\}$.  Then
    $\lie g$ is simple Lie algebra $V \cong \lie g$ and $\Kg\cong\mathbb
    R$.
  \item Assume that $\left(V{\otimes}\gp\right)^G \not= \{0\} \not= \Kg$.
    Then $(\lie g,V)$ is either $(\su(3),\mathbb{C}^3)$ or $(\lie
    g_2,V^7)$. In both cases $\Kg$ is an irreducible representation not
    isomorphic to $V$ nor to $\mathbb R$.
  \end{enumerate}
  In all cases the invariant tensors are three-forms.
  \qed
\end{lemma}

Lemma~\ref{lem:1} allows us to make the following classification Theorem:

\begin{theorem}\label{thm:3}
  Let $(M,g,V)$ be a $G$-irreducible manifold determined by a metric
  connection $\Nh$ on $(M,g)$. Assume that the torsion of $\Nh$ is parallel
  with respect to $\Nh$. Then one of the following statements holds:
  \begin{enumerate}[\upshape (a)]
  \item $(M,g)$ is locally isometric to a non-symmetric, isotropy
    irreducible homogeneous space,
  \item $(M,g)$ is locally isometric to one of the irreducible symmetric
    spaces $(G\times G)/G$ or $G^{\mathbb C}/G$,
  \item $(M,g)$ has weak holonomy $\SU(3)$ or $G_2$,
  \item the torsion of $\Nh$ vanishes and $(M,g,V)$ is the $G$-manifold
    determined by the Levi-Civita connection and $G$ is the Riemannian
    holonomy group of $(M,g)$.
  \end{enumerate}
\end{theorem}

\begin{proof}
  Assume that $(M,g,V)$ is a $G$-irreducible Riemannian manifold determined
  by a metric connection such that the torsion $\eta$ is non-trivial and
  parallel: $\Nh\eta=0$.
  
  If the space of algebraic curvature tensors is trivial then
  proposition~\ref{fn:1} applies $(M,g,\Nh)$ is an Ambrose-Singer manifold
  and locally isometric to an isotropy irreducible space.
  
  If $\Kg\not=\{0\}$ then $\lie g$ acts irreducibly on $V$ by
  Lemma~\ref{lem:5}. The torsion $\eta$ is therefore skew-symmetric by
  Lemma~\ref{lem:1}. Since $\Nh\eta=0$ we may write
  \begin{equation*}
    \RC=\Rh+(\eta^2)
  \end{equation*}
  where
  \begin{equation*}
    (\eta^2)_{X,Y}Z=[\eta_X,\eta_Y]Z-\eta_{\eta_XY-\eta_YX}Z.
  \end{equation*}
  Note that since $\eta$ is skew-symmetric we have
  \begin{equation*}
    (\eta^2)=\eta^2+\Bone\eta^2
  \end{equation*}
  where $\eta^2_{X,Y}Z := \eta_Z(\eta_XY)$. Also note that
  \begin{equation*}
    (\Bone\eta^2)_{X,Y}Z = \eta_X(\eta_YZ) - \eta_{\eta_XY}Z - \eta_Y(\eta_XZ)
    = (\eta_X.\eta)_YZ 
  \end{equation*} 
  where $.$ denotes the standard action of $\so(n)$ on $V {\otimes}
  \so(n)$. 
  
  Assume that $\Bone\eta^2 = 0 = \eta.\eta = \Bone\Rh$. Then both $\Rh$ and
  $\eta^2$ are algebraic curvature tensors. Furthermore, since $0 =
  \Bone\eta^2 = \eta.\eta$ the torsion tensor takes values in $V{\otimes}\lie
  g'$ where $\tilde{\lie g} = \stab\eta \geqslant \lie g$.  This means that
  $(\eta^2)=\eta^2\in\curv(\tilde{\lie g})$ and therefore $\eta^2$ spans a
  trivial submodule $\curv(\tilde{\lie g})$. Now Lemma~\ref{lem:1} and
  Table~\ref{tab:hol} shows that $\tilde{\lie g}$ must be a simple Lie
  algebra and that $V\cong \tilde{\lie g}$, and thus also that $\lie g =
  \tilde{\lie g}$ and $\Kg\cong\mathbb R$. So $\Rh=\kappa\eta^2$ for some
  function $\kappa\colon M\to\mathbb R$. But $(M,g)$ is Einstein by
  Theorem~\ref{thm:7}, so the scalar curvature $s_g$ is constant. But
  \begin{equation*}
    s_g = \sum_{i,j}g(\RC_{e_i,e_j}e_j,e_i) =
    (1+\kappa)\sum_{i,j}g(\eta_{e_i}e_j,\eta_{e_i}e_j) = 
    2(1+\kappa)\norm{\eta}^2
  \end{equation*}
  where $\{e_i:i=1,\dots,n\}$ is an orthonormal basis of $V$ and both $s_g$
  and $\norm{\eta}^2$ are constants. Therefore $\kappa$ must be constant
  too, whence $\LC\RC = (1+\kappa)(\Nh(\eta^2)+\eta.(\eta^2)) = 0$.

  Finally, if $\eta.\eta\not=0$ then the projection $\xi$ of $\eta$ to
  $V{\otimes}\tilde{\lie g}$ is non-trivial. Applying Lemma~\ref{lem:1}
  shows that $\tilde{\lie g}=\lie g$, $\xi=\eta$ and that $(\lie g,V)$
  is either $(\su(3),\mathbb C^3)$ or $(\lie g_2,V^7)$.
\end{proof}

\section{Invariant Intrinsic Torsion}
\label{sec:6}

Empirical evidence suggests that condition~(\ref{cond:B}) of
Theorem~\ref{thm:1} in fact implies that the intrinsic torsion must be
invariant, i.e., $W \cong k\mathbb R$. Let this serve as motivation for
considering that case in particular detail. Theorem~\ref{thm:1} with
$W\cong k\mathbb R$ becomes
\begin{proposition}
  Let $(M^n,g,V)$ be a Riemannian $G$-manifold.  Assume that the intrinsic
  torsion takes its values in the $G$-submodule $W\cong k\mathbb R\subset
  V{\otimes}\gp$.  Then $g$ is Einstein if the following conditions are
  satisfied:
  \begin{enumerate}[\upshape (a$'$)]
  \item $(V{\otimes} S^2_0V)^G=\{0\}$\label{cond:A'}
  \item $V$ is irreducible.
  \item $\lie g$ is a proper subalgebra of $\so(n)$ and $\lie
    g\not=\un(n/2)$.
  \end{enumerate}\qed
\end{proposition}

As before, these assumptions imply that the intrinsic torsion is a
three-form. 

\subsection{Invariant versus Parallel Torsion}

The assumption of invariance of the intrinsic torsion appears more general
than that of parallel torsion. It is clear that a manifold with parallel
intrinsic torsion must have invariant intrinsic torsion, since the
$G$-structure given by the holonomy of the minimal connection in this case
leaves the torsion invariant. In this section, we will prove that the
converse holds under a quite weak condition.

\begin{theorem}\label{thm:5}
  Let $(M,g,V)$ be a $G$-manifold with skew-symmetric intrinsic
  torsion taking values in some submodule $W$ of $V{\otimes}\gp$. If
  $(V{\otimes} W{\otimes} S^2W)^G=\{0\}$ then the intrinsic torsion is
  parallel with respect to the minimal connection.
\end{theorem}
\begin{proof}
  Let $\Nt$ be the minimal connection and $\xi$ the intrinsic torsion of
  $M$. Write $\gp$ for the orthogonal complement of $\lie g$ in $\so(n)$.
  We write the Riemannian curvature tensor as
  \begin{equation*}
    \RC=\Rt+(\Nt\xi)+(\xi^2).
  \end{equation*}
  where $\Rt \in \Lambda^2V{\otimes}\lie g$, $(\Nt\xi) \in
  \Lambda^2V{\otimes}\gp$ and $(\xi^2) \in (S^2(\gp))^G$. We have that
  $(\Nt\xi)$ is the anti-symmetrisation on the first two factors of $\Nt\xi
  \in V{\otimes}\Lambda^3V$. When restricted to $V{\otimes}\Lambda^3V$,
  this anti-symmetrisation is an isomorphism $V{\otimes}\Lambda^3V\cong
  \Lambda^4V+\Lambda^2(\Lambda^2V)$. Therefore,
  \begin{equation*}
    0 = \left<(\Nt\xi),\RC\right> = \left<(\Nt\xi),(\Nt\xi)+(\xi^2)\right>
    = \left<(\Nt\xi),(\Nt\xi)\right>.
  \end{equation*} 
  The last equality holds since the assumption $(V{\otimes} W{\otimes}
  S^2W)^G=\{0\}$ implies $(\Nt\xi)\in V\otimes W$ and $(\xi^2)\in S^2W$ are
  orthogonal.
\end{proof}

\begin{corollary}\label{cor:4}
  Let $(M,g,V)$ be a $G$-manifold with invariant intrinsic torsion such
  that $V^G=\{0\}$ and such that the intrinsic torsion is skew-symmetric.
  Then the intrinsic torsion is parallel with respect to the minimal
  connection.\qed
\end{corollary}

The conclusion, that $\Nt\xi=0$, implies that $H := \Hol(\Nt) \subset
\Lie{Stab}(\xi)$. However, as $\Nt$ is a $G$-connection we also have $H\subset
G$. Thus, if $\xi\not=0$ then the holonomy group $H$ must be a proper
subgroup of $G$. 

\begin{lemma}\label{lem:8}
  Let $(M,g,V)$ be a $G$-manifold with intrinsic torsion $\xi$
  and minimal connection $\Nt$. Write $H$ for the holonomy group of
  $\Nt$ and let $(M,g,V')$ be the $H$-manifold determined by the
  holonomy of $\Nt$. Then the intrinsic torsion of $(M,g,V')$ is $\xi$.
\end{lemma}
\begin{proof}
  Since $\Nt$ is a $G$-connection $H=\Hol(\Nt)\subset G$. Therefore
  $\Nt$ is a $H$-connection and, moreover, $\gp\subset\lie h^\perp$,
  whence $\xi\in V'{\otimes}\lie h^{\perp}$.  
\end{proof}

\begin{proposition}
  Let $(M,g)$ be a Riemannian manifold. Assume there exists a
  $G$-structure on $M$ with tangent representation $V$ and
  skew-symmetric intrinsic torsion taking values in $W$, where $V$ and
  $W$ satisfy $(V{\otimes} W{\otimes} S^2W)^G=\{0\}$. Then there
  exists an $H$-structure on $M$ with invariant skew-symmetric
  intrinsic torsion which is parallel with respect to the minimal
  connection.
\end{proposition}
\begin{proof}
  Assume that $(M,g,V)$ is a Riemannian $G$-manifold with intrinsic
  torsion $\xi\in
  W\subset\Lambda^3V\cap\left(V{\otimes}\gp\right)$. If
  $\left(V{\otimes}W{\otimes}S^2W\right)^G=\{0\}$ then $\Nt\xi=0$ by
  Theorem~\ref{thm:5}. As we have argued above, $H=\Hol(\Nt)$ then
  leaves $\xi$ invariant. The Proposition now follows from
  Lemma~\ref{lem:8}.
\end{proof}

\subsection{Classification}

As corollaries of Lemma~\ref{lem:1}, Theorem~\ref{thm:3} and
Corollary~\ref{cor:4} we obtain the following classifications.
\begin{theorem}\label{thm:8}
  Let $(M,g,V)$ be a $G$-irreducible Riemannian manifold. If the
  intrinsic torsion of $M$ as a $G$-manifold is invariant
  skew-symmetric and non-vanishing then either $(M,g)$ has weak
  holonomy $\SU(3)$ or $G_2$ or $(M,g)$ is locally isometric to a
  non-symmetric isotropy irreducible homogeneous space. In particular,
  $(M,g)$ is Einstein.\qed
\end{theorem}

Recall that a $G$-manifold $(M,g,V)$ with intrinsic torsion taking values
in $W\subset V\otimes\gp$ for which the three
conditions~(\ref{cond:A}),~(\ref{cond:B}) and~(\ref{cond:C}) of
Theorem~\ref{thm:1} are satisfied is $G$-irreducible and has skew-symmetric
intrinsic torsion. For the particular case of $W\cong k\mathbb R$ we now
have:

\begin{theorem}
  Let $(M,g,V)$ be a Riemannian $G$-manifold.  Assume that the intrinsic
  torsion takes its values in the $G$-submodule $W\cong k\mathbb R\subset
  V{\otimes}\gp$. Furthermore, assume that $V$ satisfies
  \begin{enumerate}[\upshape (a$'$)]
  \item $(V{\otimes} S^2_0V)^G=\{0\}$
  \item $V$ is irreducible.
  \item $\lie g$ is a proper subalgebra of $\so(n)$ and $\lie
    g\not=\un(n/2)$.
  \end{enumerate}
  Then either $(M,g)$ has weak holonomy $\SU(3)$ or~$G_2$ or $(M,g)$
  is locally isometric to a non-symmetric isotropy irreducible
  homogeneous space.\qed
\end{theorem}

\providecommand{\bysame}{\leavevmode\hbox to3em{\hrulefill}\thinspace}

\end{document}